\documentclass{amsart}
\usepackage{amsmath, amsfonts, amssymb, amstext, amscd,amsthm}

\newtheorem{sscount}{Theorem}
\newtheorem{Gmbound}[sscount]{Theorem}
\newtheorem*{sscount2}{Theorem \ref{sscount}}
\newtheorem{thm}{Theorem}[section]
\newtheorem{prop}{Proposition}[section]
\newtheorem{lem}{Lemma}[section]
\newtheorem{cor}{Corollary}[section]

\DeclareMathOperator{\Proj}{{Proj}}

\makeatletter
\def\imod#1{\allowbreak\mkern10mu({\operator@font mod}\,\,#1)}
\makeatother

\begin{document}


\title[Counting Generating Invariants]{Counting Generating Invariants \\ Under Semisimple Group and Torus Actions}
\author{Harlan Kadish}
\date{June 14, 2011}
\thanks{Supported by DMS-0502170, Enhancing the Mathematical Sciences Workforce in the 21st Century (EMSW21) Research Training Group (RTG):  Enhancing the Research Workforce in Algebraic Geometry and Its Boundaries in the Twenty-First Century}
\email{hmkadish@umich.edu}

\begin{abstract}
Fix a semisimple linear algebraic group, choose an irreducible representation of highest weight $\lambda$, and consider the irreducible representations of highest weight $n\lambda$.  As $n$ goes to infinity, the cardinality of a minimal set of generating invariants grows faster than any polynomial in $n$.  On the other hand, combinatorial methods yield sub-exponential upper bounds for the growth of generating sets for torus invariants on the binary forms.
\end{abstract}

\maketitle


\section{Introduction}

Let $k$ be an algebraically closed field and $G$ a linear algebraic group over $k$.  For a rational representation $V$ of $G$, an element $\sigma \in G$ acts on $f(x) \in k[V]$ by left translation: $\sigma \cdot f(x) = f(\sigma^{-1} x)$.  Working over $k = \mathbb C$, Hilbert in 1890 showed that the invariant subring $k[V]^G$ = $\{f \in k[V] \mid \sigma\cdot f = f\}$ is finitely generated for classical groups $G$ \cite{Hilb1890}, and in 1893 he outlined an algorithm to produce generators for $k[V]^G$ \cite{Hilb1893}.  It is now well-known that $k[V]^G$ is finitely generated for any reducitve algebraic group $G$ over an algebraically closed field $k$.  Modern computational invariant theory studies the generators of $k[V]^G$ and writes algorithms to compute them.  For example, in 1993 Sturmfels fleshed out Hilbert's algorithm for $GL_n$ \cite[ch.~4]{SturmAlgors}.  In 1999, Derksen provided an algorithm to compute generators of $k[V]^G$ for linearly reductive $G$, requiring a single Gr\"obner basis calculation \cite{DerkAlg}.  Kemper's 2003 algoritm for reductive $G$ in any characteristic requires additional normalization and integral closure procedures \cite{KemperComp}.  There also exist degree bounds for generators of $k[V]^G$ for semisimple groups, tori, and linearly reductive groups, by Popov \cite{PopovBnd}, Wehlau \cite{WehToriBnds}, and Derksen \cite{PolyBounds}, respectively.  The latter two bounds are polynomial functions of $\dim V$.

Little is known about the cardinality of minimal generating sets for invariant rings.  Estimates of such would provide bounds for the runtime of algorithms that compute invariants.  For a semisimple, linear algebraic group $G$, let $V_\lambda$ be the representation of $G$ with highest weight $\lambda$.  Parameterize with integers $n\geq 0$ the family of representations $V_{n\lambda}$ with highest weight $n\lambda$.  Let $S^d(V_{n\lambda})^G$ denote the degree-$d$ invariant polynomials on $V_{n\lambda}$.  We fix $d$ and apply a ring structure to $\oplus_{n\geq 0} S^d(V_{n\lambda})$, graded now by $n$.  It turns out that $\dim S^d(V_{n\lambda})^G$ grows like a polynomial in $n$, whose degree is a linear function of $d$.  Choosing high enough $d$, we obtain the following:

\begin{sscount}\label{sscount}
Let $\rho: G \to GL(V_\lambda)$ be a non-trivial, rational representation of highest weight $\lambda$.  The minimal cardinality of a generating set for $k[V_{n\lambda}]^G$ grows faster than any polynomial in $n$, and hence faster than any polynomial in $\dim V_{n\lambda}$.
\end{sscount}

On the other hand, we present a sub-exponential upper bound on the size of minimal generating sets for invariants of the multiplicative group $\mathbb G_m \cong k^*$.  In a representation of $\mathbb G_m$, the invariants are monomials of weight zero that are \emph{indecomposable}: not products of two other monomials of weight zero. As an application, an upper bound for the number of indecomposable invariants would provide an upper bound for a minimal generating set for the monoid of solutions $(x_1,\ldots, x_{n-1}) \in \mathbb Z^{n-1}, x_i \geq 0$ for all $i$, to the linear congruence
\[x_1 + 2 x_2 + \cdots + (n-1) x_{n-1} \equiv 0 \imod{n}.\]

\noindent Harris and Wehlau \cite{WehHar} recount the history of this problem.

To produce a stronger upper bound for generators for both $\mathbb G_m$-invariants and the solutions to linear congruences, consider the torus invariants of binary forms, as follows.  Let $\text{char}(k)=0$, and let $V_n$ denote the space $k[x,y]_n$ of homogeneous polynomials of degree $n$ in two variables, the \emph{binary forms of degree} $n$.  The algebraic group $SL_2(k)$ acts on $f(x,y) \in V$ as follows:

\[\left(\begin{array}{cc} \alpha & \beta \\ \gamma & \delta \end{array}\right) \cdot f(x,y) = f(\alpha x + \gamma y,\beta x + \delta y). \]

\noindent A maximal torus of $SL_2$ is isomorphic to $\mathbb G_m$.  Let $T$ be the maximal torus of diagonal matrices in $SL_2$. Writing a binary form as $a_n x^n + a_{n-1} x^{n-1}y + \cdots + a_0 y^n$, the coordinate ring on $V_n$ is

\[ k[V_n]\cong \begin{cases} k[ a_n, a_{n-2} , \ldots, a_0,  \ldots a_{-(n-2)}, a_{-n}] & n \text{ even} \\  k[ a_n, a_{n-2} , \ldots, a_1, a_{-1}  \ldots a_{-(n-2)}, a_{-n}] & n \text{ odd} \end{cases} \]
\noindent where $a_i$ has weight $i$ under the action of $T$.  Thus upper bounds for minimal generating sets of $k[V_n]^T$ for even $n$ would provide bounds for minimal generating sets of both
\begin{itemize}
\item $\mathbb G_m$-invariants in representations where $n/2$ is the largest magnitude of a weight.
\item Monoids of solutions to linear congruences modulo $n/2$.
\end{itemize}

Now, for an upper bound on the indecomposable invariant monomials in $k[V_n]^T$, one first counts, for each $n\geq 1$, the $S\subseteq \{-n, -(n-1), \ldots, n\}$ such that $\sum_{a\in S} a = 0$ and no subset of $S$ has this property (the ``subset sum problem").  Olson proved that the size of such $S$ is no more than $2\cdot 3\sqrt{n}$ \cite{Olson}.  In the context of monomials, this result leads to a degree bound, and one can conclude,

\begin{Gmbound} The cardinality of a minimal generating set for $k[V_n]^T$ is at most
\begin{enumerate}
\item[(a)] $O(ne^{6\sqrt{n}\log 2n})$ for odd $n$,
\item[(b)] $O(ne^{12\sqrt{n/2}\log n})$ for even $n$.
\end{enumerate}
\end{Gmbound}

Since $SL_2$ is semisimple, the two streams in this paper converge at the binary forms.  We begin with a summary of some results on counting $SL_2$-invariants of binary forms and counting solutions to linear congruences.  Indeed, the estimating technique in the proof of Theorem 1 mirrors the technique of Howe's estimates for $G = SL_2$ \cite{HoweInvars}.  For intuition and comparison, we first apply this technique to $k[V_n]^{SL_2}$ and then to arbitrary semisimple groups; in the latter case, the non-polynomial growth result appears to be new.  We then provide the combinatorial argument that counts torus invariants in Theorem 2.  

\section{History of Counting Invariants}

\subsection{$SL_2$ Invariants}
In a 1983 paper, Kac employs the ``Luna Slice Method'' over $k = \mathbb C$ to reduce questions of one representation to a ``better'' representation \cite{Kac}.  Let $V$ be a representation of a reductive linear algebraic group $G$.  For a point $p \in V$, let $G_p$ denote its stabilizer in $G$, and let $S_p$ be a $G$-stable complement to the tangent space of $G\cdot p$ in $V$.  If $G\cdot p$ is closed,
then the size of a minimal generating set for $k[V]^G$ is at least as large as the size of a minimal generating set for $k[S_p]^{G_p}$.

Kac considers the action of $G=SL_2(\mathbb C)$ on $V_n$, the binary forms of degree $n$.  First assume $n$ is odd, and choose $p = x^{n-1}y+xy^{n-1} \in V_n$.  Then Kac shows that $G_p \cong \mathbb Z_{n-2}$,
that the orbit $G\cdot p$ is closed, and that a generator $\sigma$ of $G_p$ acts on $k[S_p] = k[x_0,x_1, \ldots,x_{n-3}]$ by $\sigma \cdot x_i =\zeta_i x_i$.  Then a minimal generating set of $k[S_p]^{G_p}$ consists of monomials $m=x_0^{a_0}\cdots x_{n-3}^{a_{n-3}}$ such that
\[ a_1 + 2a_2 + \cdots + (n-3)a_{n-3} \equiv 0 \imod{n-2}\]

\noindent and such that $m$ is not divisible by another different invariant.  Let $p(k)$ be the partition function, and let $\phi(k)$ denote the number of numbers $1,\ldots, k-1$ relatively prime to $k$. Kac counts at least $p(n-2) + \phi(n-2) - 1$ generators for $k[S_p]^{G_p}$, by counting the partitions of $n-2$ and the monomials $x_i^{n-2}$ for every $i$ relatively prime to $n-2$.  Therefore, when $d$ is odd, this number also provides a lower bound for the size of a minimal generating set for $k[V_n]^{SL_2}$.  The analysis proceeds similarly for even $n$.  Now, Hardy and Ramanujan \cite{HardyRama} (and independently, Usplensky \cite{Uspensky}) found the asymptotic growth rate
\[ p(n) \sim \frac{1}{4\sqrt{3}n}e^{\pi\sqrt{2n/3}},\]

\noindent so these generating sets exhibit non-polynomial but sub-exponential growth.

Howe in 1987 computes more explicit estimates of the size of ``fundamental generating sets'' for $R_n=k[V_n]^{SL_2}$ \cite{HoweInvars}.  Let $\mathfrak m$ denote the maximal homogeneous ideal of $R_n$.  By the graded Nakayama lemma, a set $S \subset R_n$ generates $R_n$ if and only if the image of $S$ in $R_n/\mathfrak m$ generates $R_n/\mathfrak m$.  If $S$ has minimal size, then Howe calls $S$ a fundamental generating set.  He proves the following:

\begin{thm}
Let $\Gamma_n$ denote the number of fundamental invariants of $R_n$, and let $R_n(d)$ denote the degree-$d$ piece of $R_n$.  Then for fixed degree $d$,
\begin{enumerate}
\item[(a)] $\Gamma_n/(\dim R_n(d)) \to 1$ as $n\to\infty$.
\item[(b)] For constants $c_d$, the number of fundamental generators in degree $d$ for $d \geq 4$ is asymptotically
\[ \begin{cases} \frac{1}{2}(d!)^{-1}c_d n^{d-3}/(d-3) & \text{$nd$ even,} \\
0 & \text{$nd$ odd.}\end{cases} \]
\end{enumerate}
\end{thm}

\noindent Howe also provides formulas for the $c_d$ in terms of binomial coefficients.  By ``asymptotically'' Howe means ``the difference between the two expressions is small in comparison with either,'' when $n$ is large enough.  He concludes ``that almost all invariants of a fixed degree are eventually fundamental.''

\subsection{Torus and Cyclic Group Invariants} The following lemma relates the positive integer vector solutions to Kac's equation $\sum_{i=1}^{n-1} a_i \cdot i \equiv 0 \imod{n}$ to $\mathbb G_m$-invariants:

\begin{lem}
Let $T=k^*$ act on $x_i$ with weight $i$.  Identify $\mathbb Z_n$ with the $n$th roots of unity in $T$.  The evaluation homomorphism
\[ \textup{ev}\colon f(x_1, \ldots, x_n, x_{-n}) \mapsto f(x_1, \ldots, x_n, 1)\]

\noindent provides a $\mathbb Z_n$-equivariant isomorphism
\[ k[x_1, \ldots, x_n, x_{-n}]^T \to k[x_1, \ldots, x_n]^{\mathbb Z_n}.\]
\end{lem}

\begin{proof}
Surjectivity is clear.  For injectivity, it suffices to show that the ideal $(x_{-n}-1)\subset k[x_1, \ldots, x_n, x_{n-1}]$ contains no $T$-invariants.  Recall every $T$-invariant is a sum of invariant monomials.  If $f \in (x_{-n}-1)$ is a $T$-invariant, then half of the monomials of $f$ have nonzero weight, which is absurd.
\end{proof}

In a 2006 paper \cite{WehHar}, Harris and Wehlau consider the general problem, for integers $w_i$, of producing all solutions $A=(a_1,\ldots, a_r)\in\mathbb N^r$ to an equation
\[ w_1 x_1+ w_2 x_2 + \cdots + w_r x_r \equiv 0 \imod{n},\]

\noindent They note that finding solutions to this equation is equivalent to finding solutions to Kac's equation,
\[ x_1 + 2x_2 + \cdots + (n-1)x_{n-1} \equiv 0 \imod{n}.\]

\noindent To state their result, they call a solution \emph{decomposable} if, in the monoid of solutions, it can be written as a sum of two non-trivial solutions; they call it \emph{indecomposable} otherwise.  There are only finitely many indecomposable solutions: if, say, $a_i\geq n$, then one may subtract off the extremal solution $(0,\ldots,n,\ldots, 0)$ that is non-zero in the $i$th place.

The degree of a solution $A$ is $\deg(A) = \sum a_i$.  The indecomposable solutions $A=(a_1,\ldots, a_{n-1})$ correspond to generators $x_1^{a_1}\cdots x_{n-1}^{a_{n-1}}$ for $k[x_1,\ldots, x_n]^{\mathbb Z_n}$ in the proposition above.  By studying a faithful permutation action on solutions in high degree, Wehlau and Harris conclude that if $k\geq \left\lceil n/2\right\rceil+1$, then there are exactly $p(n-k)\phi(n)$ indecomposable solutions in degree $k$.  Note that $p(n-k)$ is the number of partitions of $n$ into $k$ parts. This count provides a lower bound for the number of indecomposable solutions to $\sum_{i=1}^n i\cdot x_i \equiv 0\imod{n}$.


\section{Generating Invariants for the Binary Forms}

\subsection{The Geometry of $SL_2$ Orbits} Let $k$ be an algebraically closed field, let $V_n$ denote the binary forms of degree $n$, and assume for this section that $\text{char}(k)=0$.  The results of Kac \cite{Kac} and Howe \cite{HoweInvars} already indicate that the size of a minimal generating set for $k[V_n]^{SL_2}$ grows faster than any polynomial in $n$.  Nevertheless, this example provides intuition and motivation for the study of representations of general semisimple groups.  In both cases, we first parametrize by integers $n$ a family of representations, such that the dimension of the $n$th representation is polynomial in $n$.  We then study how the number of generating invariants in high degree grows as a function of $n$.  Also in both cases, we derive our formulas from geometry:

\begin{lem} For the binary forms $V_n$ with $n \geq 3$, the generic $SL_2$ orbit is closed, of dimension 3.
\end{lem}

\begin{proof}
The non-vanishing of the discriminant gives a dense open set of forms with no double roots.  Claim the orbit of such a form is closed.  First consider the diagonal torus $T$ in $SL_2$.  The $T$-weight spaces of $V_n$ are spanned by monomials of the form $x^i y^{n-i}$.  Since $f$ has only single roots and degree at least 3, it involves monomials of both positive and negative weight. Hence if $\gamma\colon k^* \to T$ is a 1-parameter subgroup, then
\[ \lim_{t\to 0} \gamma(t)\cdot f \text{ does not exist.}\]

\noindent Now let $\gamma\colon k^* \to SL_2$ be any 1-parameter subgroup. Find $\sigma \in SL_2$ such that $\sigma \gamma \sigma^{-1}$ lies in the diagonal torus $T$ in $SL_2$.  Noting that $\sigma\cdot f$ also has all single roots,
\[ \lim_{t\to 0} \gamma(t)\cdot f = \lim_{t\to 0} \gamma(t) \sigma^{-1} \cdot \sigma f
= \sigma^{-1} \cdot \lim_{t\to 0} \sigma \gamma(t) \sigma^{-1} \cdot \sigma f, \]

\noindent which also does not exist.  By the Hilbert-Mumford Criterion, the orbit of $f$ is closed.

For $d\geq 3$, we may consider any three factors of $f$ as a triple of points in $\mathbb P^1$.  From the analysis of the complex plane, an element $\sigma \in SL_2$ is uniquely determined by its action on three distinct points, which it sends to a triple of distinct points.  Thus the stabilizer of $f$ is finite, and $\dim\overline{SL_2\cdot f} = 3$.
\end{proof}

\begin{lem} For $n \geq 3$, the categorical quotient has $\dim V_n/\!\!/SL_2 = n-2$. \end{lem}
\begin{proof}
Let $\pi\colon V_n \to V_n/\!\!/SL_2$ be the categorical quotient, a surjection of irreducible varieties.  Since the generic orbit is closed of dimension 3,
\[ 3 = \dim V_n  - \dim V_n/\!\!/SL_2 = n+1 - \dim V_n/\!\!/SL_2. \quad \]
\end{proof}


\subsection{Counting $SL_2$ Invariants} Let $V = V_1 = \{a x + b y\mid a,b \in k\}$ be the binary forms of degree 1 over an algebraically closed field $k$.  Then the space of binary forms of degree $d$ is isomorphic to $S^d(V)$, and $S^e(S^d(V))$ is isomorphic to the space of degree-$e$ regular functions on $V_d$.  That is, $S^e(S^d(V)) = k[V_d]_e$.

\begin{prop} For $V = V_1$ and natural numbers $d,e$,
\[ S^e(S^d(V)) \cong S^d(S^e(V)).\]
\end{prop}

\begin{proof}
The linear factorization of $f \in V_d$ yields a surjective, $SL_2$-equivariant morphism of varieties $\pi\colon V^d \twoheadrightarrow V_d.$  Let $S_d$ act on $V_d$ by permuting the factors, and let the torus $(k^*)^{d-1}$ act as follows:
\[ (t_1,\ldots, t_{d-1})\cdot (f_1,\ldots, f_d) = (t_1 f_1, t_1^{-1}t_2f_2, \ldots, t_{d-2}^{-1}t_{d-1}f_{d-1}, t_{d-1}^{-1}f_d).\]

\noindent Then $\pi^{-1}(f)$ is a $T\rtimes S_d$-orbit, and we have an isomorphism \[ \pi^*\colon k[V_d]\overset\sim\to k[V^d]^{T \rtimes S_d} = \left[ \overbrace{S(V) \otimes \cdots \otimes S(V)}^d \right]^{T\rtimes S_d}. \]

\noindent where $S(V)$ is the symmetric algebra.  For details, please see Derksen and Kemper's book \cite[p.~164]{DK}.

For $(f_1, \ldots, f_d) \in V^d$, write $f_i = (a_ix + b_iy)$.  If $c_0 x^d + c_1 x^dy + \cdots + c_dy^d \in V_d$, then for every $j$, $\pi^*(c_j) \in k[V^d]$ is homogeneous of degree $d$ in the $a_i, b_i$.  Since $V \cong V^*$, the map $\pi^*$ gives
\begin{eqnarray*}
S^e(S^d(V)) &\overset\sim\to& \left[ \bigoplus_{\sum e_i = de} S^{e_1}(V)\otimes \cdots \otimes S(V)^{e_d}\right]^{T\rtimes S_d} \\
&=& \left[\overbrace{S^e(V)\oplus \cdots \oplus S^e(V)}^d\right]^{S_d} \\
&=& S^d(S^e(V)).
\end{eqnarray*}
\end{proof}

\begin{lem}\label{lem:Campbell}
Let $R$ be a graded, Cohen-Macaulay domain of dimension $n$. If $R$ has Hilbert-Poincar\'e series $\sum a(d)t^d$, then there is a constant $c > 0$ such that
\[ \limsup_d\left\{\frac{a(d)}{d^{n-1}}\right\}  = c > 0,\]

\noindent and $a(d)/d^{n-1} \geq c$ for a sequence of integers $d$ with constant difference. 
\end{lem}

\begin{proof}
Let $\ell$ be the least common multiple of the degrees of a set of generators for $R$.  Let $R[\ell; i] = \oplus_m R_{m\ell + i}$, the ring of elements of degree congruent to $i$ modulo $\ell$.  Then from Section 4 of Campbell, et al. \! \cite{HilbFunct}, if $R$ is Cohen-Macaulay, then each nontrivial $R[\ell;i]$ has Hilbert polynomial $H_i(m)$ of degree $n - 1$.  What is more, if $R$ is a domain, then the leading coefficient $c$ of each nontrivial $H_i(m)$ is equal to that of $H_0(m)$; the constant $c$ is the degree of the $R[\ell; i]$.  Thus there exist infinitely many $d$, with period at most $\ell$, such that $a(d)/d^{n-1} = c+O(d^{-1})$, and the result follows.
\end{proof}

Since $SL_2$ acts linearly on $V_n$, one can find generating sets for $k[V_n]^{SL_2}$ such that each polynomial is homogeneous.  Call a subset $\Gamma$ of a $k$-algebra $R$ \emph{minimal} of it has minimal cardinality among all generating sets.  By the graded Nakayama lemma, every minimal, homogenous generating set has the same cardinality.  Recall that the method below mirrors that of Howe \cite{HoweInvars} to produce asymptopic formulas for the growth of generating sets.  Here we stop at showing that the growth is faster than polynomial, which is the extent of our result for general semsimple groups.

\begin{prop} As $n\to\infty$, the size of a minimal set of generators for $k[V_n]^{SL_2}$ grows faster than any polynomial in $n$.
\end{prop}

\begin{proof}
Suppose $\Gamma_n$ is a generating set  for $R_n := k[V_n]^{SL_2}$ of minimal cardinality.  Let $\Gamma_n(d)$ be the number of generators in $\Gamma_n$ of degree $d$.  If $k[V_n]^{SL_2}$ has Poincar\'e series $\sum a_n(d)t^d$, then for large $d$,
\[ \Gamma_n(d) \geq a_n(d) - \sum_{i=1}^{\lfloor d/2\rfloor} a_n(i)\cdot a_n(d-i) =  a_d(n)-\sum_{i=1}^{\lfloor d/2\rfloor} a_i(n)\cdot a_{d-i}(n). \]

\noindent That is, we then subtract from $\dim R_d$ the number of products of elements from $R_{<d}$, assuming no relations among them.  The result is a lower bound for the number of generators in degree $d$.  We next recall $S^i(V_n) = S^n(V_i)$ to substitute $a_n(i) = a_i(n)$.

Consider the sum on the far right above.  Now, $a_i(1) = a_1(i) = 0$ for all $n$.  When $i = 2$,  $a_i(2)=a_2(i)$ is 1 or 0, as $k[V_2]^{SL_2}$ is generated by the discriminant.  Recall $\dim R_n = n-2$.  Thus by Lemma \ref{lem:Campbell},  $\limsup_n\{ a_2(n)a_{d-2}(n)/n^{d-5}\}$ is a constant (albeit a function of $d$).  Thus $a_2(n)a_{d-2}(n)=O(n^{d-5})$.  Similarly, for $3 \leq i \leq \lfloor d/2\rfloor$, $a_i(n)\cdot a_{d-i}(n)=O(n^{i-3} \cdot n^{d-i-3}) = O(n^{d-6})$.  Thus in the relation
\[ \Gamma_n(d) \geq a_d(n)-\sum_{i=1}^{\lfloor d/2\rfloor} a_i(n)\cdot a_{d-i}(n), \]

\noindent  if $d\geq 6$, then the right-most term grows as $O(n^{d-5})$  By Lemma \ref{lem:Campbell}, there exists $c >0$ such that $a_d(n)/n^{d-3} \geq c$ for a sequence of integers $n$ with constant difference.  For $n$ in this sequence, $a_d(n) = O(n^{d-3})$ for large $n$.  Choosing $d$ arbitrarily large forces $\Gamma_n(d)$ to grow faster than any polynomial in $n$.
\end{proof}


\section{Generating Invariants of Semisimple Groups}

Let $G$ be a semisimple algebraic group over a field $k$ of characteristic 0.  Fix a Borel subgroup $B = T \ltimes U$, where $T$ is a maximal torus and $U$ is the maximal unipotent subgroup in $B$.  Let $V_\lambda$ be the representation of $G$ of highest weight $\lambda$ with respect to $T$, which is unique up to isomorphism.  We will show that when $G$ is semisimple, the cardinality of a minimal generating set of $k[V_{n\lambda}]^G$ grows faster than any polynomial in $n$.  We first describe the family of representations $V_{n\lambda}$ parametrized by $n$.  We then develop the geometry underlying our growth estimate for $\dim S^d(V_{n\lambda})$.  This estimate plays the role of the reciprocity $S^d(V_e) \cong S^e(V_d)$ for binary forms, and we can prove our first main theorem similarly.

\subsection{The Parameterization by Highest Weight}
When $U$ acts on reductive $G$ on the right, then $k[G]^U = \oplus_{\lambda\geq 0} V_\lambda$ as graded rings, where the latter is the direct sum of the irreducible representations $V_\lambda$ whose highest weight $\lambda$ is positive; for background on this ``ring of covariants,'' see Derksen and Kemper's book \cite[p.~156]{DK}.  Choose a positive weight $\lambda$, and consider the subring $R_\lambda := \oplus_{n\geq 0} V_{n\lambda}$.

\begin{lem}\label{lem:kVlammbdafg}
The ring $R_\lambda$ is finitely generated, namely, if $v_{\lambda^*}$ is the lowest weight vector of $(V_\lambda)^*$, then $R_\lambda \cong k[\overline{G\cdot v_{\lambda^*}}]$.
\end{lem}

\begin{proof}
Let $v_{\lambda^*} \in( V_\lambda)^*$ be a \emph{lowest} weight vector, of weight $-\lambda$, of the dual space to $V_\lambda$.  Claim $R_\lambda = k[\overline{G\cdot v_{\lambda^*}}]$.  Let $f$ be the image of $(v_{\lambda^*})^*$ in $k[\overline{G\cdot v_{\lambda^*}}]$.  For $n\geq 0$, the function $f^n$ has weight $n\lambda$ under $T$ and generates a $G$-module isomorphic to $V_{n\lambda}$ in $k[\overline{G\cdot v_{\lambda^*}}]$, whence $R_\lambda \hookrightarrow k[\overline{G\cdot v_{\lambda^*}}]$.

To obtain the reverse inclusion, consider the orbit map $G \to \overline{G\cdot v_{\lambda^*}}$ defined by $g \mapsto g\cdot v_{\lambda^*}$.  Because this map is dominant, it gives rise to a $G$-equivariant injection $k[\overline{G\cdot v_{\lambda^*}}]\hookrightarrow k[G]$.  Consider the stabilizer in $T$ of $v_{\lambda^*}$,
\[ T_{\lambda^*} = \{ t \in T \mid t\cdot v_{\lambda^*} = v_{\lambda^*}\}.\]

\noindent If $h_\mu \in k[\overline{G\cdot v_{\lambda^*}}]$ is a highest weight vector of weight $\mu$, claim $\mu(T_{\lambda^*})=\{1\}$.  First note that if $h_\mu(v_{\lambda^*})=0$, then $h_\mu(T\cdot v_{\lambda^*})  = T\cdot h_\mu(v_{\lambda^*})= \{0\}$, because $h_\mu$ is a weight vector.  Let $U^-$ be the opposite unipotent subgroup to $U$ with respect to $T$; then $v_{\lambda^*}$ is $U^-$ invariant, because $U^-$ lowers the weights of $T$.  It follows
\[ \{0\} = h_\mu(v_{\lambda^*}) = h_\mu(T\cdot v_{\lambda^*}) = h_\mu(U T U^-\cdot v_{\lambda^*})\]

\noindent because $h_\mu$ is $U$-invariant as a highest weight vector.  Since $U T U^-$ is dense in $G$, $h_\mu$ would be identically zero on $G\cdot v_{\lambda^*}$, which is absurd.  Thus $h_\mu(v_{\lambda^*})\neq 0$, and for $t \in T_{\lambda^*}$,
\[ \mu(t) h_\mu(v_{\lambda^*})= t\cdot h_\mu(v_{\lambda^*}) = h_\mu(t^{-1}\cdot v_{\lambda^*}) =h_\mu(v_{\lambda^*}).\]

\noindent Therefore $\mu(T_{\lambda^*})=\{1\}$ and $\mu = n\lambda$, so every irreducible $G$-submodule of $k[\overline{G\cdot v_{\lambda^*}}]$ is one of the $V_{n\lambda}$.
\end{proof}

\noindent For comparison to the size of a generating set, note that the dimension of $V_{n\lambda}$ grows as a polynomial in $n$:

\begin{lem}\label{lem:WeylFormula}
Let $r$ be the number of positive roots of a reductive group $G$.  Then $\dim V_{n\lambda} = O(n^r)$.
\end{lem}

\begin{proof}
Let $\Phi$ be the set of roots of $G$, $\delta  = \frac{1}{2}\sum_{\alpha\succ 0}\alpha$ a sum over the positive roots, and $(\cdot, \cdot)$ an inner product on the space spanned by $\Phi$, preserved by the Weyl group of reflections.  Then by Weyl's formula in (see \cite[p.~139]{HumphLA}),
\[ \dim V_{n\lambda} = \frac{\prod_{\alpha\succ 0}(n\lambda + \delta, \alpha)}{\prod_{\alpha\succ 0} (\delta,\alpha)}.\]

\noindent The number $r$ of positive roots of $G$ satisfies $2r + \dim T = \dim G$.
\end{proof}


\subsection{Generic Closed Orbits in Cartesian Products} For a finite-dimensional vector space $V$ over an algebraically closed field $k$, let $\rho\colon G \to GL(V)$ be a non-trivial, rational representation of the semisimple algebraic group $G$.

\begin{lem}\label{lem:GinSLV}
If $G$ is semisimple and $\rho: G \to GL(V)$ is a representation, then the image of $\rho$ lies in $SL(V)$.  If $\rho$ is not trivial, then $\dim\rho(G) \geq 2$.
\end{lem}

\begin{proof} The image $\rho(G)$ of $G$ in $GL(V)$ is also semisimple, and
\[\rho(G) = [\rho(G),\rho(G)] \subseteq [GL(V),GL(V)] = SL(V).\]

\noindent If $\rho$ is non-trivial, then $\dim\rho(G) \geq 1$, but there are no connected, semisimple algebraic groups of dimension 1 (see \cite[p.~131]{Humph}).
\end{proof}

\begin{lem}\label{lem:genxdorbit}
For an $n$-dimensional vector space $V$, let $X = \mathbb P(V)$.  Let $d \geq n+1$ and let $SL(V)=SL_n$ act diagonally on $X^d$.  If $\widetilde{X^d}$ is the affine cone over $X^d$, then the generic orbit of $SL_n$ acting on $\widetilde{X^d}$ is closed, of dimension $\dim SL_n$.
\end{lem}

\begin{proof}
Fixing a basis for $V$, let $f\colon V^{d} \to k$ be the product of the $(n\times n)$-minors of an $n\times d$ matrix.  This $f$ defines a function on $X^d$ and also on the affine cone $\widetilde{X^d}$.  Choose $p \in \widetilde{X^d}$ with $f(p)\neq 0$, and let $[p]$ be its image in $X^d$.  Then $[p]$ defines $d$ points in $\mathbb P(V)$, no $n$ of which lie in the same hyperplane.  Thus the stabilizer of $[p]$ in $SL_n$ is finite (namely, the scalar matrices of $SL_n$ such that the product of the entires is 1).  It follows that the stabilizer of $p$ is finite, whence $\dim (SL_n\cdot p) = n^2-1.$ This dimension holds for the orbit of generic $p$ with $f(p)\neq 0$.

Note that $f$ is an invariant function on $\widetilde{X^d}$, because the $SL(V)$ fixes determinants.  So if $q \in \widetilde{X^d}$ lies in the orbit closure of $p$, then $f(q) = f(p)$.  Thus $\dim (SL_n\cdot q) = n^2-1$ as well.  Since orbits in the boundary of $SL_n\cdot p$ must have strictly smaller dimension, the point $q$ must lie in the orbit of $p$.  Therefore, the generic orbit is closed.
\end{proof}

Recall that $V_\lambda$ is a highest-weight representation of a semisimple group $G$, and $R = \oplus_{n\geq 0}V_{n\lambda}$.  Let $Z = \Proj R$, and consider the sum of tensors over $k$,
\[ C = \bigoplus_{n\geq 0} \overbrace{V_{n\lambda} \otimes \cdots \otimes V_{n\lambda}}^{d}. \]

\noindent The $d$th Cartesian product of $Z$ is
\[ Z^d = \Proj C = \overbrace{Z \times_k \cdots \times_k Z}^{d}. \]

\noindent As $G \subseteq \text{Aut}(R)$, $G$ acts rationally on $Z$, hence diagonally on $Z^d$, hence on the affine cone $\widetilde{Z^d}$.

\begin{lem}\label{lem:genzdorbit}
In the above notation, if $\rho: G \to GL(V_\lambda)$ is an irreducible representation of highest weight $\lambda$ and $d > \dim V_\lambda$, then the generic orbit of $G$ acting on $\widetilde{Z^d}$ is closed, of dimension $\dim \rho(G)$.
\end{lem}

\begin{proof}
Suppose $\dim V_\lambda = n$.  Note $Z$ is a subvariety of $X = \mathbb P(V_\lambda^*)$: indeed, if $k[V_\lambda] = k[x_1,\ldots, x_n]$, then there is a surjection $k[x_1, \ldots, x_n] \twoheadrightarrow \oplus_{m\geq 0}V_{m\lambda}$ by sending the $x_i$ onto an ($n$-dimensional) basis for $V_\lambda$.  Thus $\widetilde{Z^d}$ is a closed subvariety of $\widetilde{X^d}$, and $\widetilde{Z}$ spans $V_\lambda^*$ because $k[\widetilde{Z}]$ contains $V_\lambda$.

As in the proof of Lemma \ref{lem:genxdorbit}, let $f\colon V_\lambda^{d} \to k$ be the product of the $(n\times n)$-minors of an $n\times d$ matrix.  Since $\widetilde{Z}$ is irreducible and spans $V_\lambda$, the generic $p \in \widetilde{Z^d} \subseteq \widetilde{X^d}$ has $f(p) \neq 0$.  Then for generic $p \in \widetilde{Z^d}$, the orbit $SL(V_\lambda)\cdot p$ is closed in $\widetilde{X^d}$.  Recall that the stabilizer $SL(V_\lambda)_p$ of $p$ is a finite set of scalar matrices, whence normal. It follows that the orbit $SL(V_\lambda)\cdot p$ is isomorphic to $SL(V_\lambda)/SL(V_\lambda)_p$ as a variety, so it is an algebraic group.  By Lemma \ref{lem:GinSLV}, $\rho(G)$ is a closed subgroup of $SL(V_\lambda)$ (see \cite[p.~54]{Humph}), and the $\rho(G)$ action on $V_\lambda$ factors through the $SL(V_\lambda)$ action.  Thus for generic $p \in \widetilde{Z^d}$, $G\cdot p$ is isomorphic to $\rho(G)/\rho(G)_p$, the homomorphic image of an algebraic group. Thus $G\cdot p$ is closed in $SL(V_\lambda)\cdot p$, with dimension $\dim\rho(G)$.  Since $G\cdot p \subseteq \widetilde{Z^d}$, the result follows.
\end{proof}


\subsection{Counting Generating Invariants}

To count generating invariants for large $n$, we again need to understand the degree $d$ component of $k[V_{n\lambda}]^G$.

\begin{prop}
Let $\rho: G\to GL(V_\lambda)$ be a non-trivial, rational representation of highest weight $\lambda$.  Write $S^d(V_{n\lambda})^G \cong k[V_{n\lambda}]_{d}^G$, the degree-$d$ homogeneous piece of $k[V_{n\lambda}]^G$.  Let $m = \dim \rho(G)$.  Then there is an integer $c$ with $1 \leq c < m$ such that for large $n$,
\[ \dim S^d(V_{n\lambda})^G \leq O(n^{cd-1}),\]

\noindent and when both $n$ and $d$ are large,
\[ \dim S^d(V_{n\lambda})^G = O(n^{cd-m}),\]
\end{prop}

\begin{proof}
Letting the symmetric group $S_d$ permute the $d$ factors of each $n$-graded piece of
\[ \bigoplus_{n\geq 0} \overbrace{V_{n\lambda} \otimes \cdots \otimes V_{n\lambda}}^{d}, \]

\noindent take the categorical quotient,
\[ Z^d /\!\!/ S_d = \Proj\left(\bigoplus_{n\geq 0} S^d(V_{n\lambda})\right).\]

\noindent Next take the quotient by the $G$ action on each copy of $X$:
\[ \left(Z^d /\!\!/ S_d\right)/\!\!/G = \Proj\left(\bigoplus_{n\geq 0} S^d(V_{n\lambda})^G\right). \]

\noindent These actions of $G$ and $S_d$ commute.  From Lemma \ref{lem:genzdorbit}, if $d \geq \dim V_\lambda+1$, then the generic orbit of $G$ acting on the cone $\widetilde{Z^d}$ is closed and of dimension $\rho(G)=m$.  Therefore, for large enough $d$,
\[ \dim(Z^d/\!\!/S_d)/\!\!/G = d\cdot \dim Z - m.\]

\noindent Now, Lemma \ref{lem:kVlammbdafg} yields that
\[ \dim Z + 1= \dim R = \dim k[\overline{G\cdot v_{\lambda^*}}] \leq \dim \rho(G) = m.\]

\noindent Let $c = \dim Z$.  Note $c \geq 1$, because for large $d$, $\widetilde{Z^d}$ contains an orbit of dimension $m > 1$.  Thus the Hilbert polynomial for $\oplus_n S^d(V_{n\lambda})^G$ has degree $cd - m$ for large $d$, and degree bounded by $cd-1$ otherwise.
\end{proof}

As above, let $k[V_{n\lambda}]^G$ have Hilbert-Poincar\'e series $\sum_{d=0}^\infty a_n(d)t^d$.

\begin{sscount2}
Let $\rho: G \to GL(V_\lambda)$ be a non-trivial, rational representation of highest weight $\lambda$.  The minimal cardinality of a generating set for $k[V_{n\lambda}]^G$ grows faster than any polynomial in $n$, and hence faster than any polynomial in $\dim V_{n\lambda}$.
\end{sscount2}

\begin{proof}
Let $\Gamma_n$ denote the minimal cardinality of a generating set of $k[V_{n\lambda}]^G$, and let $N = \dim V_\lambda$.  From the proof above, if $d > N$, then $\dim S^d(V_{n\lambda})^G = O(n^{cd-m})$ for large $n$ and a constant $c$ with $1\leq c < m$.  Then for large $n$,
\begin{eqnarray*}
\Gamma_n &\geq& a_n(d) - \sum_{i=1}^{\lfloor d/2 \rfloor} a_n(i)a_n(d-i) \\
&=& a_n(d) - \sum_{i=1}^{N} a_n(i)a_n(d-i) - \sum_{i=N+1}^{\lfloor d/2 \rfloor} a_n(i)a_n(d-i) \\
&\approx & a_n(d) -\sum_{i=1}^{N} a_n(i)a_n(d-i) - \sum_{i=N+1}^{\lfloor d/2 \rfloor} n^{ci-m}n^{c(d-i)-m} \\
\end{eqnarray*}

\noindent where the approximation symbol indicates an asymptotic estimate for sufficiently large $n$ and $d$. When $1\leq i \leq N$ and $n$ is large, we bound $a_n(i) = O(n^{ci - 1})$.  Assume $d - N > N, m$, so that for $i \leq N$ we may bound $a_n(d-i) = O(n^{c(d-i)-m})$.  Then for such $i \leq N$, we have $a_n(i)a_n(d-i) = O(n^{cd-m-1})$, and we obtain that for large $n$ and $d$,
\begin{eqnarray*}
O(\Gamma_n) &\geq& a_n(d) - N\cdot n^{cd - m - 1} - d n^{cd-2m} \\
&\approx& n^{cd-m}-N\cdot n^{cd - m - 1} - d n^{cd-2m} \\
&\approx& n^{cd-m} \\
\end{eqnarray*}

\noindent Fixing $d$ arbitrarily large, it follows that the size of a minimal generating set for $k[V_{n\lambda}]^G$ grows faster than any polynomial in $n$.  The final assertion of the theorem follows because, by Lemma \ref{lem:WeylFormula}, $\dim V_{n\lambda}$ grows like a polynomial in $n$.
\end{proof}


\section{An Upper Bound for $\mathbb G_m$ Invariants}

Recall that torus invariants are generated by monomials.  We say an invariant monomial is \emph{indecomposable} if it is not the product of non-constant invariants.  Consider the action of $T\cong k^*$ on a polynomial ring $k[x_1, \ldots, x_n, x_{-n}]$ defined by
\[ t\cdot x_i = t^i x_i.\]

\noindent Then the invariants are monomials of weight zero, whose positive-weight part is a multiple of $n$.

\begin{prop}
The cardinality of a minimal generating set for
\[ k[x_1, \ldots, x_n, x_{-n}]^T \quad \text{is} \quad O\left(e^{6\sqrt{n}\log 2n}\right).\]
\end{prop}

\begin{proof}
A minimal generating set contains only monomials of the form $x_nx_{-n}$ and $mx_{-n}^k$, where $m \in (x_1, \ldots, x_{n-1})$ is properly divisible by no monomial of weight congruent to $0$ modulo $n$.  Now, Olson shows that if $S \subset \mathbb Z_n$ has order at least $3\sqrt{n}$, then a subset of $S$ has trivial sum \cite{Olson}.  Now, a minimal generating set of $k[x_1, \ldots, x_n, x_{-n}]^T$ can be chosen such that each (monomial) generator properly includes no invariant.  So such a set can be chosen such that each generating monomial includes no more than $3\sqrt{n}$ distinct variables,though possibily with repetition.

An algorithm of Derksen and Kemper to construct torus invariants implies a degree bound of $2n-1$ for a generating set of $k[x_1, \ldots, x_n,x_{-n}]^T$, by computing within the convex hull of the variables' weights \cite[p.~159]{DK}
.  This linear bound may only hold when the torus has rank 1; a more general bound appears in Wehlau \cite{WehToriBnds}. An upper bound for the number of generators in degree $d$ is
\[ \binom{n}{\lfloor 3\sqrt{n}\rfloor} \cdot \binom{\lfloor 3\sqrt{n}\rfloor + d-1}{\lfloor 3\sqrt{n}\rfloor -1} \leq n^{\lfloor 3\sqrt{n}\rfloor} \cdot (\lfloor 3\sqrt{n}\rfloor + d -1)^{\lfloor 3\sqrt{n}\rfloor-1}. \]

\noindent The first term on the left counts ways of choosing $3\sqrt{n}$ variables; the second term counts monomials of degree $d$ with $3\sqrt{n}$ variables.  Summing the upper bound over degrees $d$ up to $2n-1$ yields
\[ O\left((n^{3\sqrt{n}} \cdot (3\sqrt{n} + 2n -2)^{3\sqrt{n}}\right) \leq O((2n)^{6\sqrt{n}}) = O(e^{6\sqrt{n}\log 2n}).\]
\end{proof}

%

\noindent Recall our initial interest in $k[V_n]^T$, the the torus-invariant functions on the binary forms of degree $n$.  Note that $k[V_{2n}]^T$ is isomorphic to
\[ B_n := k[x_{-n},x_{-n+1}, \ldots, x_0, \ldots, x_n]^T.\]

\begin{prop} The cardinality of a minimal generating set for $B_n$ is
\[ O\left(ne^{12\sqrt{n}\log 2n}\right).\]
\end{prop}

\begin{proof}
For $r \geq 1$, let $mx_{-r}^k$ be an invariant monomial such that
\[m \in k[x_{r-1}, \ldots, x_1, x_{-1}, \ldots, x_{-(r-1)}]\]

\noindent and $m$ is properly divisible by no monomial of weight congruent to $0$ modulo $r$.  Since $x_i$ and $x_{-r+i}$ have the same weight modulo $r$, Olson's theorem yields that $m$ involves no more than $2\cdot 3\sqrt{r}$ distinct variables, which may occur with multiplicity.  Otherwise, if $m'$ divides $m$ and $m'$ involves $6\sqrt{r}$ distinct variables, then $m'$ is divisible by a monomial of weight congruent to $0$ modulo $r$.

Recall that monomials of the form $mx_{-r}^k$ in a minimal generating set have degree at most $2n-1$, for any $r$.  Since an upper bound for the number of generating $mx_{-r}^k$ in degree $d$ is
\[ \binom{2r}{6\sqrt{r}} \cdot \binom{6\sqrt{r} + d-1}{6\sqrt{r}-1} \leq (2r)^{6\sqrt{r}} \cdot (6\sqrt{r} + d -1)^{6\sqrt{r}-1}, \]

\noindent summing these upper bounds up to degree $2n-1$ yields
\[ O\left((2r)^{6\sqrt{r}} \cdot (6\sqrt{r} + 2n -2)^{6\sqrt{r}}\right) \leq O((2n)^{12\sqrt{r}}) = O(e^{12\sqrt{r}\log 2n}).\]

\noindent Repeat the argument for invariants of the form $x_r^km$ with
\[ m \in k[x_{r-1}, \ldots, x_1, x_{-1}, \ldots, x_{-(r-1)}].\]
In either case, the weight of $m$ determines the exponent $k$.  Note that the number of invariants $x_rx_{-r}$ grows linearly, and these monomials together generate the invariants. The result follows by choosing the largest $r =n$.
\end{proof}

\noindent In the invariant $mx_r^k$ considered in the above proof, the monomial $m$ may be divisible by invariants not involving $x_r$.  Nevertheless, the upper bound holds for monomials whose highest-weight variable is $x_r$, and the following corollary is our Theorem 2:

\begin{cor} The cardinality of a minimal generating set for $k[V_n]^T$ is
\begin{enumerate}
\item[(a)] $O(ne^{6\sqrt{n}\log 2n})$ for odd $n$,
\item[(b)] $O(ne^{12\sqrt{n/2}\log n})$ for even $n$.
\end{enumerate}
\end{cor}

\begin{proof}  For odd $n$ there is a $T$-equivariant isomorphism
\[ k[V_n]^T \cong k[x_{-n}, x_{-n+2}, \ldots, x_{-1}, x_1, \ldots, x_n]^T,\]

\noindent where $x_i$ has torus weight $i$.  Following the argument above, the invariant $mx_r^k$, say, involves no more than $3\sqrt{r}$ variables among the $r+1$ variables $x_{r-2}, \ldots, x_{-(r-2)}$ of distinct weight modulo $r$.  Making these adjustments to the above calculations, but retaining the degree bound $2n-1$ (from the convex hull of the variables' weights), yields an upper bound for the size of a minimal generating set:
\[ (r+1)^{3\sqrt{r}} \cdot (3\sqrt{r} + 2n -2)^{3\sqrt{r}} \leq O((2n)^{6\sqrt{r}}) = O(e^{6\sqrt{r}\log 2n}).\]

\noindent The weight $r$ varies from $1,3, 5 \ldots, n$, and the result follows for odd $n$,.

For even $n$, the isomorphism
\[ k[V_n]^T \cong k[x_{-n/2}, x_{-n/2+1}, \ldots, x_0, \ldots, x_{n/2}]^T,\]

\noindent where $x_i$ has torus weight $i$, makes way for the previous proposition.
\end{proof}

\bibliographystyle{plain}
\bibliography{thesis}   

\end{document}